\documentclass{amsart}
\usepackage{amsmath, amsthm, amssymb, amscd}

\DeclareMathOperator{\as}{asdim} \DeclareMathOperator{\Ind}{Ind}
\DeclareMathOperator{\asI}{asInd}
\DeclareMathOperator{\diam}{diam} 
\DeclareMathOperator{\dist}{dist}

\DeclareMathOperator{\mesh}{mesh}

\begin{document}
\author[G.C. Bell]{G. C. Bell}
\address{Department of Mathematics, Penn State University,
University Park, PA 16802} \email{bell@math.psu.edu}

\author[A.N. Dranishnikov]{A. N. Dranishnikov}
\address{Department of
Mathematics, PO Box 118105, Gainesville, FL 32611-8105}
\email{dranish@math.ufl.edu}

\title[A Hurewicz-type theorem for asymptotic dimension]{A
Hurewicz-type theorem for asymptotic dimension and applications to
geometric group theory}

\newcommand{\To}{\longrightarrow}
\newcommand{\sA}{\mathcal{A}}
\newcommand{\sB}{\mathcal{B}}
\newcommand{\sU}{\mathcal{U}}
\newcommand{\sV}{\mathcal{V}}
\newcommand{\sW}{\mathcal{W}}
\newcommand{\sH}{\mathcal{H}}
\newcommand{\sY}{\mathcal{Y}}
\newcommand{\sZ}{\mathcal{Z}}
\newcommand{\Real}{\mathbf{R}}
\newcommand{\Complex}{\mathbf{C}}
\newcommand{\so}{\Rightarrow}
\newcommand{\sX}{\mathcal{X}}
\newtheorem{Theorem}{Theorem}
\newtheorem*{thm}{Theorem}
\newtheorem*{cor}{Corollary}
\newtheorem*{prop}{Proposition}
\newtheorem{Lemma}{Lemma}
\newtheorem{Proposition}{Proposition}
\newtheorem{Corollary}[Theorem]{Corollary}
\newtheorem*{defn}{Definition}
\bibliographystyle{amsplain}
\keywords{asymptotic dimension, free products, nilpotent groups}

\thanks{2000 {\em Mathematics Subject Classification.} Primary:
20F69, 20F65; Secondary: 20E08, 20E06\\
The second author was partially supported by NSF Grant
DMS-0305152}

\begin{abstract}  We prove an asymptotic analog of the classical
Hurewicz theorem on mappings which lower dimension.  This theorem
allows us to find sharp upper bound estimates for the asymptotic
dimension of groups acting on finite dimensional metric spaces and
allows us to prove a useful extension theorem for asymptotic
dimension. As applications we find upper bound estimates for the
asymptotic dimension of nilpotent and polycyclic groups in terms
of their Hirsch length.  We are also able to improve the known
upper bounds on the asymptotic dimension of fundamental groups of
complexes of groups, amalgamated free products and the
hyperbolization of metric spaces possessing the Higson property.
\end{abstract}

\maketitle

\section{Introduction}  In classical dimension theory the
Hurewicz theorem on mappings which lower dimension is a powerful
tool.  One statement of the theorem is the following (see
\cite{En95}, for example).

\begin{thm} Let $X$ and $Y$ be compact metric spaces and $f:X\to Y$ a
continuous map.  Suppose that there is some n so that for every
$y\in Y,$ $\dim f^{-1}(y)\le n.$  Then $\dim X\le \dim Y+n.$
\end{thm}

Gromov defined the asymptotic dimension of a metric space in his
study of asymptotic invariants of finitely generated groups in
\cite{Gr93}. The {\em asymptotic dimension} of a metric space X,
$\as X,$ is defined to be the smallest integer $n$ so that for
every $R$ there is a uniformly bounded cover of $X$ so that no
$R$-ball in $X$ meets more than $n+1$ elements of the cover.

Asymptotic dimension is a coarse invariant, see \cite{Ro3}, so in
particular it is a quasi-isometry invariant.  Much of the interest
in asymptotic dimension has been directed at showing that certain
classes of groups have finite asymptotic dimension in a natural
metric, see for example: \cite{Be04}, \cite{BD1}, \cite{DJ},
\cite{Gr93}, \cite{Ji}. For the most part we are concerned with
finitely generated groups in this note. The metric we associate to
a finitely generated group is the word metric coming from some
finite generating set. As any two such metrics are Lipschitz
equivalent, $\as \Gamma$ is well-defined without reference to a
generating set.  Finite $\as$ results for groups became important
following a theorem of Yu \cite{Yu1} which showed that the Novikov
higher signature conjecture holds for $\Gamma$ with $\as
\Gamma<\infty.$ Because of Yu's theorem many results were aimed
only at showing that $\as \Gamma$ was finite and not concerned
with computing the exact dimension.

There are analogies between local topology and asymptotic topology
as well as ``dictionaries" which seek to translate between the
local and asymptotic worlds, see \cite{Dr00}. For dimension,
however, the correspondence is often not direct. For example, in
classical dimension theory we have the Urysohn-Menger Theorem
which gives the sharp upper bound: $\dim(X\cup Y)\le \dim X+\dim
Y+1;$ in the asymptotic case, by contrast, the authors showed in
\cite{BD1} (also, see section 4) that $\as(X\cup Y)= \max\{\as
X,\as Y\}.$  Also, $\cup_{i=1}^\infty \dim F_i= \max\{\dim F_i\}$
for closed sets, but not every finitely generated group $\Gamma$
has $\as\Gamma=0,$ so there cannot be a direct analog of the
countable union theorem.  In section 4 we state a countable union
theorem for $\as$ from \cite{BD1}.

In contrast with the union theorems, the asymptotic analog of the
Hurewicz theorem on mappings which lower dimension is very much a
direct analog of its classical counterpart. In the third section
we prove our version of an asymptotic Hurewicz theorem:

\noindent{\bf Theorem \ref{Hurewicz}.} {\em Let $f:X\to Y$ be a
Lipschitz map of a geodesic metric space to a metric space.
Suppose that for every $R>0,$ $\{f^{-1}(B_R(y))\}_{y\in Y}$
satisfies the inequality $\as \le n$ uniformly (see section 2).
Then $\as X\le \as Y+ n.$}

This leads to an upper bound estimate for finitely generated
groups acting on finite dimensional metric spaces which agrees
with the formula conjectured in \cite{BD1} as well as \cite{Ro3}.

The fourth section is devoted to a specialized version of the
Hurewicz theorem.  When the codomain of the Lipschitz map is a
tree and some asymptotic disjointness conditions are satisfied,
the estimate on the asymptotic dimension can be improved.  This
leads to a generalization of the formula for the asymptotic
dimension of a free product of groups.  In particular, we consider
a free product of pointed metric spaces and compute their
asymptotic dimension.  This also leads to an upper bound estimate
on the asymptotic dimension of free products with amalgamation
which, unfortunately, is given in terms of the asymptotic
dimension of quotients.  Whereas $\as$ is monotonic in subsets,
the $\as$ of a quotient cannot be determined from the $\as$ of the
spaces in the quotient, (see the remarks following Theorem
\ref{amalgam}).

In the final section we apply the Hurewicz theorem to finitely
generated groups acting by isometries on metric spaces.  This
leads us to an extension theorem for $\as$:

\noindent{\bf Theorem \ref{extension}.}  {\em Let $\phi:G\to H$ be
a surjection of a finitely generated group with $\ker\phi=K.$
Then, $\as G\le \as H+\as K.$}

We conclude the paper with applications of our extension theorem
to nilpotent groups, some amalgamated free products and the
hyperbolization of metric spaces possessing the Higson property.

\section{Asymptotic Dimension and Uniform Mapping Cylinders}

As mentioned in the introduction, Gromov defined the $\as$ of a
metric space in \cite{Gr93}.  There are many equivalent
formulations of asdim.  As we will need to pass between them, we
summarize some of the equivalences below.  (All of these
equivalences are stated in \cite{Gr93}, for explicit proofs, see
\cite{BD2}, \cite{Ro3}.)

\begin{thm} For a metric space $X$
the following are equivalent:
\begin{enumerate}
    \item for every $D>0$ there exist $D$-disjoint families
    $\sU_0,\ldots,\sU_n$ of uniformly bounded sets whose union covers $X;$
    \item for every $R>0$ there exists a uniformly bounded
    cover $\sU$ of $X$ in
    which no $R$-ball in $X$ meets more than $n+1$ elements of the
    cover $\sU$;
    \item for every $L>0$ there exists a uniformly bounded
    cover of $X$ with multiplicity $\le n+1$ and with Lebesgue
    number $>L$;
    \item for every $\epsilon>0$ there is a uniformly
    cobounded, $\epsilon$-Lipschitz map $\phi:X\to K$ to a
    uniform polyhedron of dimension $\le n.$
\end{enumerate}
\end{thm}

Recall that a family $\sU$ of subsets of a metric space is said to
be $D$-{\em disjoint} if $d(U,V)>D$ for every $U\neq V$ in $\sU.$
The condition on the cover in item (2) of the theorem is often
referred to as {\em R-multiplicity} $\le n+1.$   The {\em Lebesgue
number} of a cover $\sU$ of a metric space $X$ is
$L(\sU)=\inf\{\max\{d(x,X\setminus U)\mid U\in\sU\}\mid x\in X\}$.
A {\em uniform polyhedron} is the geometric realization of a
simplicial complex in $\ell^2,$ with the metric it inherits as a
subset.  Finally, a map $\phi$ to a uniform polyhedron is {\em
uniformly cobounded} if there is a number $B$ so that
$\diam(\phi^{-1}(\sigma))\le B$ for all simplices $\sigma.$

\begin{defn} A metric space $X$ has $\as X\le n$ if it satisfies
any of the equivalent conditions of the previous theorem.
\end{defn}

Often we will need to work with the canonical projection of a
cover to its nerve. In fact, the implication $(3)\implies (4)$ of
the previous theorem can be seen by simply applying the canonical
projection to the nerve.  Let $\sU$ be an open cover of a metric
space $X$. The {\em canonical projection to the nerve} $p:X\to
Nerve(\sU)$ is defined by the partition of unity
$\{\phi_U:X\to\mathbb{R}\}_{U\in\sU}$, where
$\phi_U(x)=d(x,X\setminus U)/\sum_{V\in\sU} d(x,X\setminus V)$.
The family $\{\phi_U:X\to\mathbb{R}\}_{U\in\sU}$ defines a map $p$
to the Hilbert space $\ell^2(\sU)$ with basis indexed by $\sU$.
The nerve $Nerve(\sU)$ of the cover $\sU$ is realized in
$\ell^2(\sU)$ by taking every vertex $U$ to the corresponding
element of the basis. Clearly, the image of $p$ lies in the nerve.

In \cite{BD2} the authors showed that the canonical projection
$p:\sU\to Nerve(\sU)$ of a cover $\sU$ with multiplicity $k+1$ and
Lebesgue number $L$ is $\frac{(2k+3)^2}{L}$-Lipschitz, (cf.
Proposition \ref{nu}).

In the statement of the asymptotic Hurewicz theorem we need the
following natural notion of uniformity for $\as$ defined by the
authors in \cite{BD1}. A family $\{X_\alpha\}$ of subsets of a
metric space $X$ satisfies the inequality $\as X_\alpha\le n$ {\em
uniformly} (see \cite{BD2}) if for large $D>0$ there is an $R>0$
such that there exist $R$-bounded, $D$-disjoint families
$\sU_\alpha^0,\ldots,\sU^n_\alpha$ so that
$\cup_{i=0}^n\sU^i_\alpha$ is a cover of each $X_\alpha.$ A basic
example of families satisfying $\as X_\alpha\le n$ uniformly is
when all the families are isometric.

In view of the fact that any tree $T$ has $\as T\le 1,$ (see
\cite{Ro3}) the authors proved in \cite{BD2} what could be called
a first approximation to the asymptotic Hurewicz theorem. In
particular, the main result was

\begin{thm}\cite[Theorem 1]{BD2} Suppose that the finitely
generated group $\Gamma$ acts cocompactly by isometries on a tree
$X.$ Then, $\as \Gamma\le k+1,$ where $\as\Gamma_x\le k$ for all
stabilizers $\Gamma_x$ of vertices $x\in X.$
\end{thm}

The proof used the characterization of asymptotic dimension in
terms of uniformly cobounded, Lipschitz maps to uniform polyhedra.
The argument here is similar and relies heavily on the notion of
simplicial mapping cylinders.  We summarize the pertinent results
on simplicial mapping cylinders in the following proposition.

\begin{Proposition} \label{Prop3} \cite[Propositions 2, 3]{BD2}  For
every simplicial map $f:X\to Y$ from a $n$-dimensional simplicial
complex X the mapping cylinder $M_f$ admits a triangulation with
the set of vertices equal to the disjoint union of vertices of $X$
and $Y;$ there is a constant $c_n$ so that the quotient map
$q:X\times [0,1]\to M_f$ is $c_n$-Lipschitz, where $M_f$ is given
the uniform metric it inherits from $\ell^2.$
\end{Proposition}

\begin{Proposition}  \label{Prop4} \cite[Proposition 4]{BD2}
Let $A\subset W\subset X$ be subsets in a geodesic metric space
$X$ such that the $r$-neighborhood $N_r(A)$ is contained in $W$
and let $f:W\to Y$ be a continuous map to a metric space $Y.$
Assume that the restrictions $f|_{N_r(A)}$ and $f|_{W\setminus
N_r(A)}$ are $\epsilon$-Lipschitz.  Then $f$ is
$\epsilon$-Lipschitz.
\end{Proposition}

We end this section with a computation we will need later.

\begin{Proposition}  \label{product} Let $X,Y$ and $Z$ be metric spaces.  Suppose
$f:X\to Y$ and $g:X\to Z$ are Lipschitz functions with Lipschitz
constants $\lambda_f$ and $\lambda_g,$ respectively.  Then, the
map $f\times g:X\to Y\times Z$ defined by $x\mapsto (f(x),g(x))$
is $\sqrt{2}\max\{\lambda_f,\lambda_g\}$-Lipschitz in the product
metric: $\sqrt{d_Y^2+d_Z^2}$
\end{Proposition}

\begin{proof}  The proof is an elementary calculation:
\[\begin{array}{rcl}
d_{Y\times Z}\Big(\big(f(x),g(x)\big),\big(f(x'),g(x')\big)\Big)
& =   &\sqrt{[d_Y(f(x),f(x'))]^2+[d_Z(g(x),g(x'))]^2}\\
& \le &\sqrt{(\lambda_f^2+\lambda_g^2)d_X(x,x')^2}\\
& \le &\sqrt{2}\max\{\lambda_f,\lambda_g\}d_X(x,x')
\end{array}\]
\end{proof}

\section{An Asymptotic Hurewicz Theorem}

We need a version of Lemma 1 from \cite{BD2} (also see our Lemma
\ref{L1}). The result is very technical, so we break it up over
the next two lemmas and one proposition.

\begin{Lemma} Let $f:X\to Y$ be a $\lambda$-Lipschitz map of a
geodesic metric space to a metric space with $\lambda\ge 1.$  Let
$r>1$ be given and suppose that $\sW$ is a uniformly bounded cover
of $Y$ with uniformly bounded $\lambda r$-multiplicity.  Let
$\tau$ be a simplex in $Nerve(N_{\lambda r}(\sW))$ maximal with
respect to containment, and take $\tau'$ to be a simplex in
$\beta^1\tau$ with $d=\dim\tau=\dim\tau'.$  For $i=0,\ldots, k,$
let $W_i$ denote the vertex of $\tau'$ corresponding to an
$i$-face of $\tau.$ Put $X_{\tau'}=f^{-1}(\cup_{i=0}^k W_i).$
Finally suppose that there exist families $\sU_0,\ldots, \sU_k$ of
uniformly bounded sets with multiplicities $\le n+1$ such that
\begin{enumerate}
    \item $\sU_i$ covers $f^{-1}(\cup_{j=0}^i W_i)$ and
    \item for all $i<j$ there exist simplicial maps
$\psi^{(j)}_{(i)}:Nerve(\sU_i)\to Nerve(\sU_j).$
\end{enumerate}
Then there exists a uniformly cobounded, Lipschitz map
$\phi:X_{\tau'}\to K_{\tau'}$ to a uniform polyhedron of dimension
$n+k.$
\end{Lemma}

{\bf Remark:} The existence of such a map is not difficult to see
through use of the finite union theorem from \cite{BD1} and the
fourth definition of $\as$ given above, but we will need specific
properties of the map we construct here.

\begin{proof}  For $x\in X_{\tau'},$ define
\[t_i(x)=\max\{0,\dfrac{\lambda
r-\dist(f(x),W_i)}{\lambda r}\}.\]  Observe that $0\le t_i(x)\le
1$ with $t_i(x)=1$ precisely when $f(x)\in W_i,$ and $t_i(x)>0$ if
and only if $f(x)$ is in the interior of $N_{\lambda r}(W_i).$
Also, observe that on $X_{\tau'},$ $t_0\equiv 1.$

We define the map $\phi:X_{\tau'}\to K_{\tau'}$ as a combination
of simpler maps. First, we define $\phi_0:X_{\tau'}\to
Nerve(\sU_0)$ by $\phi_0(x)=p_{\sU_0}(x),$ where $p_{\sU_0}$
denotes the canonical projection to the nerve $Nerve(\sU_0).$  We
define $\phi_1$ before passing to a general description of
$\phi_i.$

First, put $g_1=\psi_{(0)}^{(1)}:Nerve(\sU_0)\to Nerve(\sU_1).$
Let $M_{g_1}$ denote the uniform mapping cylinder and let $q_1:
Nerve(\sU_0)\times [0,1]\sqcup Nerve(\sU_1)\to M_{g_1}$ be the
quotient and uniformization map.

Define $\phi_1:X_{\tau'}\to M_{g_1}$ by
\[\phi_1(x)=\left\{%
\begin{array}{ll}
    q_1(\phi_0(x),2t_1(x)), & \hbox{if $t_1(x)\in[0,\frac12]$;} \\
    2(1-t_1(x))\psi_{(0)}^{(1)}\phi_0(x)+(2t_1(x)-1)p_{\sU_1}(x), & \hbox{otherwise.} \\
\end{array}%
\right.\]

More generally, suppose that $\phi_{p-1}$ and $g_{p-1}$ have been
defined.  Define $g_p:M_{g_{p-1}}\to Nerve(\sU_p)$ by
\begin{multline*}
    g_p([z,t_1,\ldots,t_{p-1}])=\\
    t_{p-1}\psi^{(p)}_{(p-1)}\psi_{(0)}^{(p-1)}(z)+(1-t_{p-1})
[t_{p-2}\psi^{(p)}_{(p-2)}\psi^{(p-2)}_{(0)}(z)+(1-t_{p-2})[\cdots]].
\end{multline*}  Here we have extended the $\psi$ by defining $\psi^{(j)}_{(0)}(z)=z$ for all $z\in
Nerve(\sU_j).$

Next, put
\[\phi_p(x)=\left\{%
\begin{array}{ll}
    q_p(\phi_{p-1}(x),2t_p(x)), & \hbox{if $t_p(x)\in[0,\frac12]$;} \\
    2(1-t_p(x))\psi_{(p-1)}^{(p)}\phi_{p-1}(x)+(2t_p(x)-1)p_{\sU_p}(x), & \hbox{otherwise,} \\
\end{array}%
\right.\] where, as before, $q_p$ is the uniformization and
quotient map to the mapping cylinder of $M_{g_p}.$  Put
$\phi:X_{\tau'}\to K_{\tau'}$ equal to $\phi_k.$

First we show that $\phi$ is uniformly cobounded.  To this end,
let $\sigma\in K_{\tau'}$ be a simplex.  Suppose that
$\xi,\eta\in\sigma$ and that $x_\xi\mapsto \xi,$ $x_\eta\mapsto
\eta$ under $\phi.$  Then, $p_{\sU_k}(x_\xi)$ and
$p_{\sU_k}(x_\eta)$ lie in the same simplex and obviously
$f(x_\xi)$ and $f(x_\eta)$ lie in the same simplex. Thus
$\dist(x_\xi,x_\eta)\le \max\{2b(\sU_i),2b(\sW)\},$ which is a
uniform bound.  (Here $b(\sU_i)$ denotes an upper bound on the
diameters of the sets in $\sU_i.$)

It remains to show that the map $\phi$ is Lipschitz and compute
the Lipschitz constant.  We consider $\phi:X_{\tau'}\to M_k.$
Observe that $x\in N_{\lambda r/2}(W_k)$ if and only if $\frac12
\le t_k(x)\le 1.$  So, applying Proposition \ref{Prop4}, we see
that $\phi_k$ is Lipschitz if it is Lipschitz when $t_k(x)\in
[0,\frac12]$ and when $t_k(x)\in [\frac12, 1].$ (Here we are using
the fact that $X$ is geodesic.) But, the definitions of these maps
depend on $\phi_{k-1},$ which in turn depend on $\phi_{k-2}.$
Thus, we begin with $\phi_0$ and work up inductively.

The map $\phi_0$ is just $p_{\sU_0}$ so by \cite[Proposition
1]{BD2} it is $\frac{(2n+3)^2}{L(\sU_0)}\:$-Lipschitz, where
$L(\sU_0)$ is the Lebesgue number of $\sU_0.$  Next, we consider
$\phi_1.$  We recall the definition:
\[\phi_1(x)=\left\{%
\begin{array}{ll}
    q_1(\phi_0(x),2t_1(x)), & \hbox{if $t_1(x)\in[0,\frac12]$;} \\
    2(1-t_1(x))\psi_{(0)}^{(1)}\phi_0(x)+(2t_1(x)-1)p_{\sU_1}(x), & \hbox{otherwise.} \\
\end{array}%
\right.\] Notice that in the second case, we have
$\dist(x,f^{-1}(W_1))\le \frac{\lambda r}2,$ so by Proposition
\ref{Prop4} it suffices to show that the map is Lipschitz in both
cases; then $\phi_1$ will be Lipschitz with constant equal to the
max of the constants from each of the cases.

In the first case, we know that $\phi_0$ is
$\frac{(2n+3)^2}{L(\sU_0)}\:$-Lipschitz, $t_1(x)$ is
$\frac{2}{\lambda r}$-Lipschitz, and $q_1$ is $c_n$-Lipschitz.
Thus, the map $\phi_1$ is
$c_n\sqrt2(\max\{\frac{(2n+3)^2}{L(\sU_0)},\frac{2}{\lambda
r}\})$-Lipschitz, by Proposition \ref{product}.  In the second
case, we apply the Leibnitz rule to see that the sum is
$\frac2{\lambda r}+2\frac{(2n+3)^2}{L(\sU_0)}+\frac2{\lambda
r}+2\frac{(2n+3)^2}{L(\sU_1)}$-Lipschitz.  Hence in this case it
has Lipschitz constant equal to $2\frac{2}{\lambda
r}+\frac{(2n+3)^2}{L(\sU_0)}+ \frac{(2n+3)^2}{L(\sU_1)}.$  So, we
conclude that $\phi_1$ is $\lambda_1$-Lipschitz, where
$\lambda_1=\max\{c_n\sqrt2(\max\{\frac{(2n+3)^2}{L(\sU_0)},\frac{2}{\lambda
r}\}),2\frac{2}{\lambda r}+2\frac{(2n+3)^2}{L(\sU_0)}+
2\frac{(2n+3)^2}{L(\sU_1)}\}.$

Similarly, assuming $\phi_{p-1}$ is $\lambda_{p-1}$-Lipschitz, we
consider
\[\phi_p(x)=\left\{%
\begin{array}{ll}
    q_p(\phi_{p-1}(x),2t_p(x)), & \hbox{if $t_p(x)\in[0,\frac12]$;} \\
    2(1-t_p(x))\psi_{(p-1)}^{(p)}\phi_{p-1}(x)+(2t_p(x)-1)p_{\sU_p}(x), & \hbox{otherwise.} \\
\end{array}%
\right.\]  As before in the top case we see that the map is
$c_{n+p-1}\sqrt2\max\{\lambda_{p-1},\frac2{\lambda
r}\}$-Lipschitz.  In the second case, we apply the Leibnitz rule
again to see that the Lipschitz constant is $2\frac{2}{\lambda
r}+2\lambda_{p-1}+ 2\frac{(2n+3)^2}{L(\sU_p)}.$  Thus, we conclude
that $\phi_p$ is $\lambda_p$-Lipschitz with
$\lambda_p=\max\{c_{n+p-1}\sqrt2\max\{\lambda_{p-1},\frac2{\lambda
r}\}, 2\frac{2}{\lambda r}+2\lambda_{p-1}+
2\frac{(2n+3)^2}{L(\sU_p)}\}.$  Thus, $\phi$ is Lipschitz, with
Lipschitz constant $\lambda_k.$

\end{proof}

\begin{Lemma} \label{welldef} In the notation of the previous lemma, suppose
$\sigma$ and $\tau$ are simplices in $N'$ both of which are
maximal with respect to containment.  Suppose that
$\sigma\cap\tau=\varrho.$  Then $\phi^{(\tau)}|_{\varrho}=
\phi^{(\sigma)}|_{\varrho}.$
\end{Lemma}

\begin{proof} Suppose that the vertices of $\sigma$ are denoted
$v_0,\ldots, v_c$ and the vertices of $\tau$ are $w_0,\ldots, w_d$
where the index of the vertex corresponds to the dimension of the
cell of which the vertex is the barycenter.

First, we show that for all $x$ which map to $\varrho,$
$t^{(\sigma)}_p(x)=t^{(\tau)}_p(x)$ for all $p.$  If
$\varrho=[v_{i_0},\ldots,v_{i_s}],$ as a subsimplex of $\sigma,$
then $\varrho=[w_{i_0},\ldots,w_{i_s}],$ as a subsimplex of $\tau$
since the indices must correspond to dimensions of cells in $N.$
Clearly $t_j(x)=0$ for all the indices which do not appear in the
description of $\varrho.$  All the other $t_{i_j}$ must agree as
they are defined in terms of distances which are intrinsic to the
simplex $\varrho.$

We prove the lemma by induction on the dimension of $\varrho.$ To
begin, suppose that $\varrho$ is a point. If $\varrho=v_p=w_p,$
then $t^{(\sigma)}_p(x)=1=t^{(\tau)}_p(x)$ for any $x$ with
$p_{N'}f(x)=\varrho.$  It is also clear that if $i\neq p,$ then
$t^{(\sigma)}_i(x)=t^{(\tau)}_i(x)=0.$ Thus,
$\phi^{(\sigma)}(x)=q_{p+1}(\phi^{(\sigma)}_p(x),0,\ldots, 0).$
But $\phi^{(\sigma)}_p(x)=p_{\sU_p}(x).$ On the other hand, when
we compute $\phi$ thinking of $\varrho$ as a subsimplex of $\tau,$
we see $\phi^{(\tau)}(x)= (p_{\sU_p}(x),0,\ldots,0).$ Thus, we
obtain $\phi^{(\sigma)}(x)\approx p_{\sU_p}(x)\approx
\phi^{(\tau)}(x),$ where $\approx$ denotes the natural
identification in the mapping cylinder.

Next, consider $\varrho$ with $\sigma$-vertices $\{v_{i_0},\ldots,
v_{i_{s}}\}.$  Since the indices on the vertices agree, the
$\tau$-vertices must be $\{w_{i_0},\ldots, w_{i_{s}}\}.$ Applying
the definition, we see that $\phi^{(\sigma)}(x)=
q_{i_{s}+1}(\phi^{(\sigma)}_{i_{s}}(x),0,\ldots, 0),$
where \begin{multline*}\phi^{(\sigma)}_{i_{s}}(x)=\\
\left\{%
\begin{array}{ll}
    q_{i_{s}}(\phi^{(\sigma)}_{i_{s}-1}(x),2t_{i_{s}}(x)), & \hbox{if $t_{i_s}(x)\in[0,\frac12]$;} \\
    2(1-t_{i_{s}}(x))\psi_{(i_{s}-1)}^{(i_{s})}\phi^{(\sigma)}_{i_{s}-1}(x)+(2t_{i_{s}}
    (x)-1)p_{\sU_{i_{s}}}(x), & \hbox{otherwise.} \\
\end{array}%
\right.\end{multline*}  Similarly, we find $\phi^{(\tau)}(x)=
q_{i_{s}+1}(\phi^{(\tau)}_{i_{s}}(x),0,\ldots, 0),$ where
\begin{multline*}\phi^{(\tau)}_{i_{s}}(x)=\\
\left\{%
\begin{array}{ll}
    q_{i_{s}}(\phi^{(\tau)}_{i_{s}-1}(x),2t_{i_{s}}(x)), & \hbox{if $t_{i_s}(x)\in[0,\frac12]$;} \\
    2(1-t_{i_{s}}(x))\psi_{(i_{s}-1)}^{(i_{s})}\phi^{(\tau)}_{i_{s}-1}(x)+(2t_{i_{s}}
    (x)-1)p_{\sU_{i_{s}}}(x), & \hbox{otherwise.} \\
\end{array}%
\right.\end{multline*} By the inductive hypothesis applied to the
simplex whose vertices are $\{v_{i_0},\ldots v_{i_{s-1}}\}$ and
$\{w_{i_0},\ldots w_{i_{s-1}}\},$ we see that the maps
$\phi_p^{(\sigma)}$ and $\phi_p^{(\tau)}$ agree for all $p<i_s,$
up to identification in the mapping cylinders. Thus, the maps
$\phi^{(\sigma)}$ and $\phi^{(\tau)}$ agree on $\varrho.$
\end{proof}

\begin{Proposition} \label{Lip} Let $\epsilon>0$ be given.  Suppose $\lambda$
is a constant, $\lambda\ge 1.$  Finally, suppose that
$r<L(\sU_0)<\cdots<L(\sU_k)$ in the notation of Lemma 1, where
$r\ge \frac1\epsilon(2n+3)^26^kc_nc_{n+1}\cdots c_{n+k-1}.$ Then,
$\phi$ is $\epsilon$-Lipschitz, where the $c_i$ are the constants
from Proposition \ref{Prop3}.
\end{Proposition}

\begin{proof}  Again the proof is a simple computation.  For
$0\le p\le k,$ we show that \[\lambda_p\le
\frac{(2n+3)^2}{r}\max\{c_n\sqrt2,6\}\cdots\max\{c_{n+p-1}\sqrt2,6\}.\]
Then, $\lambda_k=
\frac{(2n+3)^2}{r}\max\{c_n\sqrt2,6\}\cdots\max\{c_{n+k-1}\sqrt2,6\},$
and so \[\lambda_k\le\frac{\epsilon (2n+3)^2
6^{k-\ell}\sqrt{2}^\ell c_{n_1}\cdots c_{n_\ell}}{6^k c_n\cdots
c_{n+k-1}},\] for some $\ell.$ Since the maps $q$ restrict to an
isometry on $t=0,$ we have $c_i\ge 1,$ and so we see that
$\lambda_k$ does not exceed $\frac{\epsilon}{(3\sqrt2)^\ell}\le
\epsilon.$

To prove the claim, we use induction.  Let $p=1.$  We saw in the
proof of Lemma 1 that
$\lambda_1=\max\{c_{n}\sqrt2\frac{(2n+3)^2}{r},
c_n\sqrt2\frac{2}{\lambda r}, 2\frac{2}{\lambda
r}+2\frac{(2n+3)^2}{L(\sU_0)}+ 2\frac{(2n+3)^2}{L(\sU_1)}\}.$ This
does not exceed $\max\{c_{n}\sqrt2\frac{(2n+3)^2}{r},
2\frac{2}{r}+2\frac{(2n+3)^2}{r}+ 2\frac{(2n+3)^2}{r}\}\le
\max\{c_{n}\sqrt2\frac{(2n+3)^2}{r}, 6\frac{(2n+3)^2}{r}\}.$ Thus,
$\lambda_1\le\max\{c_{n}\sqrt2, 6\}\frac{(2n+3)^2}{r},$ as
desired.

To prove the inductive step, we use the estimate in Lemma 1:
\[\begin{array}{rcl}
\lambda_p&= &
\max\{c_{n+p-1}\sqrt2\max\{\lambda_{p-1},\frac2{\lambda r}\},
2\frac{2}{\lambda r}+2\lambda_{p-1}+
2\frac{(2n+3)^2}{L(\sU_p)}\}\\
& \le &
\max\{c_{n+p-1}\sqrt2\lambda_{p-1},\frac4r+2\lambda_{p-1}+\frac{2(2k+3)^2}{r}\}\\
& \le & \max\{c_{n+p-1}\sqrt2,6\} \lambda_{p-1}.
\end{array}\] The last inequality follows from observing that both
$\frac2r$ and $\frac{(2n+3)^2}{r}$ are not more than
$\lambda_{p-1}.$

\end{proof}

The next proposition is another technical result which relies
heavily on the uniform inequality $\as\le n$ for a family of
metric spaces.  We also need the notion of $d$-{\em saturated
union}. Let $\sU$ and $\sV$ be families of subsets of a metric
space $X.$ Denote by $N_d(V;\sU)$ the union of $V$ and all
$U\in\sU$ with $d(V,U)\le d.$  The $d$-saturated union,
$\sV\cup_d\sU$ is defined to be the family $\{N_d(v;\sU)\mid
V\in\sV\}\cup\{U\in\sU\mid d(U,V)>d,\text{ for all } V\in\sV\}.$

\begin{Proposition} \label{bound} Let $\{F_\alpha\}$ be a collection of
subspaces of the metric space $X$ satisfying $\as F_\alpha\le n$
uniformly. Then, for any $m\in\mathbb{Z}$ and for any $L>0$ there
is a bound $b$ so that there is a $b$-bounded cover of
$\cup_{i=1}^m F_{\alpha_i}$ with multiplicity $\le n+1$ and
Lebesgue number $\ge L.$
\end{Proposition}

\begin{proof}  First, for each $m,$ we prove that the collection of
$\{\cup_{\alpha\in I} F_\alpha\}_{|I|=m}$ has asymptotic dimension
$\le n$ uniformly.

We proceed inductively, the base case being true by assumption.
For any collection of $m$ sets, write $\cup_I F_\alpha$ as
$F_{\alpha_0}\cup \bigcup_{I'} F_\alpha$ where $I'$ is the index
set $I$ with $\alpha_0$ removed.  Take $d$-disjoint, $R$-bounded
families $\sU_0,\ldots,\sU_n$ covering $F_{\alpha_0}$ and
$5R$-bounded, $r$-bounded families $\sV_0,\ldots,\sV_n$ covering
$\cup_{I'}F_\alpha,$ by the inductive hypothesis.  Then, put
$\sW_i=\sV_i\cup_d \sU_i,$ the $d$-saturated union.  Then, this
family is $d$-disjoint and $r+2(d+R)$-bounded.  Since this
construction was independent of the indexing set, the family
$\{\cup_{\alpha\in I} F_\alpha\}_{|I|=m}$ has asymptotic dimension
$\le n$ uniformly.

Now, we prove the assertion of the proposition.  Let $L$ be given.
Take $d=2L.$  Then, construct the covers $\sW_i$ as in the
preceding paragraph.  Thus, $\sW_i$ are $d$-disjoint and
$r+2(d+R)$-bounded.  Finally, put $\sW=\cup_{i=0}^n
N_{d/2}(\sW_i).$  Then $\sW$ covers the union $\{\cup_{\alpha\in
I} F_\alpha\}_{|I|=m},$ and $\sW$ has multiplicity $\le n+1$ and
the Lebesgue number of $\sW$ is greater than $d/2=L,$ as desired.
Clearly, $b=r+2R+3d$ is a uniform bound on the diameters of the
elements of $\sW.$
\end{proof}

We are finally in a position to prove the main result of the
paper, our Hurewicz-type theorem for $\as.$

\begin{Theorem} \label{Hurewicz} Let $f:X\to Y$ be a Lipschitz map of a
geodesic metric space to a metric space.  Suppose that for every
$R>0,$ $\{f^{-1}(B_R(y))\}_{y\in Y}$ satisfies the inequality $\as
\le n$ uniformly.  Then $\as X\le \as Y+n.$
\end{Theorem}

\begin{proof}  Suppose $\as Y\le k.$ For a given $\epsilon>0$
we will construct a uniformly cobounded, $\epsilon$-Lipschitz map
$\Phi:X\to K$ to a uniform simplicial complex of dimension $n+k.$

Suppose that $f$ is $\lambda$-Lipschitz. So that we can apply
Proposition \ref{Lip}, we observe that $\lambda$ can be taken to
be at least $1.$ Take $r$ as in Proposition \ref{Lip} and let
$\sW$ be a cover of $f(X)$ by uniformly bounded sets with
multiplicity $\le k+1$ whose $\lambda r$-enlargements also have
multiplicity $\le k+1.$ (Using the first definition of asymptotic
dimension it is not difficult to see that such a cover must
exist.)

Since the $W\in\sW$ are uniformly bounded, there is an $R>0$ so
that for every $W\in\sW$ there is a $y_W\in f(X)$ so that
$N_{\lambda r}(W)\subset B_R(y_W).$  Thus, $f^{-1}(N_{\lambda
r}(W)\subset f^{-1}(B_R(y_0)).$  Since by assumption, $\as
f^{-1}(B_R(y))\le n$ uniformly, we conclude that $\as
f^{-1}(N_{\lambda r}(W))\le n.$

So, for each $W\in\sW,$ let $\sU^W_0$ denote a family of covers of
$f^{-1}(N_{\lambda r}(W))$ which are uniformly bounded, have
multiplicity $\le n+1$ and have Lebesgue number $L(\sU_0)>r.$
Inductively define covers $\sU^{\sigma}_i$ covering
$f^{-1}(\sigma)$ for all $i$-dimensional simplices $\sigma$ in
$Nerve(N_{\lambda r}(\sW)).$  Insist that the covers be uniformly
bounded, have multiplicity $\le n+1$ and take
$L(\sU_i)>b(\sU_{i-1}).$   This final condition is possible by
Proposition \ref{bound}.

Let $N'$ denote the barycentric subdivision of $Nerve(N_{\lambda
r}(\sW)).$  Let $\tau$ be a simplex in $Nerve(N_{\lambda r}(\sW))$
maximal with respect to containment.  Let $\tau'$ be a simplex in
the barycentric subdivision of $\tau$ with dimension equal to that
of $\tau.$  Then we have covers $\sU_i$ for the vertices $v_i$ of
$\tau'.$  The conditions on the Lebesgue numbers and bounds of the
covers mean that there are simplicial maps
$\psi^{(j)}_{(i)}:Nerve(\sU_i)\to Nerve(\sU_j)$ whenever $i<j.$

Apply the lemma to obtain a map $\phi:X_{\tau'}\to K_{\tau'}.$ The
map is $\epsilon$-Lipschitz by Proposition \ref{Lip} and uniformly
cobounded by Lemma 1.  Glue the $K_{\tau'}$ using the face
relation on $N'.$  By Lemma \ref{welldef}, the $\phi$ agree along
common faces.  So by Proposition \ref{Prop4} they define an
$\epsilon$-Lipschitz, uniformly cobounded map $\Phi:X\to K.$
Obviously $\dim K=n+k.$
\end{proof}

We immediately obtain a result for groups acting by isometries on
trees.  This upper bound is the product-type estimate one expects,
cf. \cite{BD1}, \cite{Ro3}.

\begin{Theorem} \label{main} Let $\Gamma$ be a finitely generated group
acting on the metric space $X$ by isometries.  Fix a point $x_0\in
X.$  Suppose that $\as X\le k$ and $\as W_R(x_0)\le n,$ for all
$R.$  Then, $\as\Gamma\le n+k.$
\end{Theorem}

\begin{proof}  Fix a finite generating set $S=S^{-1}$ for
$\Gamma.$ Then, $|\Gamma|_S$ is a geodesic metric space. Define
the map $\pi:\Gamma\to X$ by $\pi(\gamma)= \gamma x_0.$ Put
$\lambda= \max\{d_X(sx_0,x_0)\mid s\in S\}.$  We claim that the
$\pi$ is $\lambda$-Lipschitz.  Since the metric space $\Gamma$ is
discrete geodesic, it suffices to check the Lipschitz condition on
pairs of points at distance $1$ from each other.  Such a pair is
of the form $(\gamma,\gamma s),$ for some $s\in S.$  Now,
\[d_X(\pi(\gamma),\pi(\gamma s))=d_X(\gamma x_0,\gamma
sx_0)=d_X(x_0,s x_0)\le \lambda.\] Thus, $\pi$ is
$\lambda$-Lipschitz.  So that we can apply Proposition \ref{Lip},
we observe that $\lambda$ can be taken to be at least $1.$

Next, observe that $W_R(x_0)=\pi^{-1}(B_R(x_0))$ and since the
action is isometric, $\gamma B_R(x_0)=B_R(\gamma x_0).$  So,
$\pi^{-1}(B_R(\gamma x_0))=W_R(\gamma x_0)$ and $W_R(x_0)$ is
isometric to $W_R(\gamma x_0)$ for all $\gamma\in\Gamma.$  Thus,
$\as \pi^{-1}(B_R(x))\le n$ uniformly.  We apply Theorem
\ref{Hurewicz} to get $\as\Gamma\le n+k.$
\end{proof}

In \cite[Lemma 1]{Be04}, the first author proved the following
result about complexes of groups.  Complexes of groups are a
natural generalization of the Bass-Serre theory of graphs of
groups.  Whereas graphs of groups describe groups acting on trees,
in the general theory trees are replaced with higher-dimensional
analogs of trees called small categories without loops (briefly
scwols).  There is not always an action associated to a complex of
groups, so the theory is not as nice as the Bass-Serre theory.
When there is an associated action, the complex of groups is
called {\em developable.}  For more details see \cite[Chapter
III.$\mathcal{C}$]{BH}.

\begin{Proposition} Let $\pi$ be the fundamental group of a
complex of groups $G(\sY)$ where $\sY$ is finite, the local groups
$G_\sigma$ are finitely generated and $\as G_\sigma\le n$ for all
$\sigma.$  Fix some $\sigma_0\in\sY.$  Then for every $R>0,$ $\as
W_R(G_{\sigma_0})\le n.$
\end{Proposition}

Applying Theorem \ref{main} and the previous proposition we
immediately obtain the following sharpening of \cite[Theorem
3]{Be04}.

\begin{Theorem}  Let $\Gamma$ be the fundamental group of a
finite, developable complex of groups corresponding to an action
by isometries on the geometric realization of the scwol $\sX.$
Suppose the local groups are finitely generated and $\as
G_\sigma\le n.$  If $\as |\sX|\le k,$ then $\as \Gamma\le n+k.$
\end{Theorem}

When $\sX$ is one-dimensional, we recover the main theorem from
\cite{BD2}.

\section{Lipschitz Mappings to Trees}

We want to generalize the main result of \cite{BDK} which gives a
formula for the asymptotic dimension of a free product of groups.

We will need the following results:

\begin{Proposition}\cite[Proposition 1]{BD2} \label{nu} For every
$K$ and every $\epsilon>0$ there exists a number
$\nu=\nu(\epsilon,k)$ such that for every cover $\sU$ of a metric
space $X$ of order $\le k+1$ with Lebesgue number $L(\sU)>\nu$ the
canonical projection to the nerve, $p_\sU :X\to Nerve(\sU)$ is
$\epsilon$-Lipschitz.
\end{Proposition}

\begin{Lemma}\cite[Lemma 1]{BD2} \label{L1} Let $A$ be a closed
subset of the geodesic metric space $X.$  Let $r>8\epsilon$ and
let $\sV$ and $\sU$ be covers of the $r$-neighborhood $N_r(A)$ by
uniformly bounded open sets such that $\sV$ has order $\le n+1,$
$Nerve(\sV)$ is orientable, and
$L(\sU)>b(\sV)>L(\sV)\ge\nu(\epsilon/4c_{n},n),$ where $c_n$ is
the constant of uniformization from Proposition \ref{Prop3}.  Then
there is an $\epsilon$-Lipschitz map $f:N_r(A)\to M_g$ to the
mapping cylinder supplied with the uniform metric of a simplicial
map $g:Nerve(\sV)\to Nerve(\sU)$ between the nerves such that $f$
is uniformly cobounded, $f|_{\partial N_r(A)}=q(p_{\sV}|_{\partial
N_r(A)},0),$ and $f|_A=p_{\sU}|_A,$ where $p_{\sU}:N_R(A)\to
Nerve(\sU)$ and $p_{\sV}:N_r(A)\to Nerve(\sV)$ are the canonical
projections to the nerves.
\end{Lemma}

To improve the Hurewicz-type estimate for maps to trees we need
the notion of asymptotic inductive dimension.  Asymptotic
inductive dimension, $\asI$, was defined by the second author in
\cite{Dr01} in order to establish connections between $\as X$ and
$\Ind \nu X$ where $\nu X$ is the Higson corona of $X.$

Let $\varphi:X\to\mathbb{R}$ be a function defined on a metric
space $X.$  For every $r>0$ let
$V_r(x)=\sup\{|\varphi(y)-\varphi(x)|\colon y\in B_r(x)\}.$  Such
a function $\varphi$ is called {\em slowly oscillating} if for
every $r>0$ and every $\epsilon>0,$ there exists a compact set
$K\subset X$ so that $V_r(x)<\epsilon$ for all $x\in X\setminus
K.$  Let $\bar X$ be the compactification of $X$ corresponding to
the family of all continuous bounded slowly oscillating functions.
The {\em Higson corona of} $X$ is the remainder $\nu X=\bar
X\setminus X.$

For any subset $A\subset X$ denote by $A'$ the trace of $A$ on
$\nu X,$ i.e. the intersection of the $\bar X$-closure of $A$ with
$\nu X.$

Let $X$ be a proper metric space.  A subset $W\subset X$ is called
an {\em asymptotic neighborhood} of $A\subset X$ if for any $x_0,$
$\lim_{r\to \infty}d(A\setminus B_r(x_0),X\setminus W)=\infty.$

Two sets $A,B\subset X$ are {\em asymptotically disjoint} if for
any $x_0,$ $\lim_{r\to\infty}d(A\setminus B_r(x_0),B\setminus
B_r(x_0))=\infty.$  Thus, $A$ and $B$ are asymptotically disjoint
precisely when their traces $A'$ and $B'$ are disjoint.

A subset $C\subset X$ is an {\em asymptotic separator} for the
asymptotically disjoint sets $A,B\subset X$ if its trace $C'$ is a
separator for $A'$ and $B'.$  Define $\asI X=-1$ if and only if
$X$ is bounded and otherwise, define $\asI X\le n$ if for any
asymptotically disjoint $A$ and $B$ in $X$ there exists an
asymptotic separator $C$ with $\asI C\le n-1.$

\begin{thm} \cite[Theorem 3]{Dr01} Let $X$ be a proper metric space with
bounded geometry and suppose that $\as X$ is positive and finite.
Then $\asI X=\as X.$
\end{thm}

There is a small problem with bounded sets.  If $K$ is a bounded
set $\asI K=-1,$ whereas $\as K=0.$  (Notice that there are
unbounded sets, for example $\{2^n\}\subset \mathbb{Z}$, with
$\as$ zero.)  In any case, for metric spaces $X$ with bounded
geometry, we have $\asI X\le \as X,$ see \cite{DrZ}.

\begin{Theorem} \label{totrees} Let $f:X\to T$ be a Lipschitz map of the geodesic
metric space $X$ with bounded geometry to a tree.  Suppose that
for any disjoint bounded sets $W$ and $W'$ in $T,$ the sets
$f^{-1}(W)$ and $f^{-1}(W')$ are asymptotically disjoint.  Suppose
that for every $R>0,$ the family $\{B_R(v)\}_{v\in T}$ satisfies
the inequality $\as f^{-1}(B_R(v))\le n$ uniformly in $v\in T,$
with $n\ge 1.$ Then $\as X\le n.$ Moreover if $f^{-1}(B_R(v))=n$
for some $v$ and $R,$ then $\as X=n.$
\end{Theorem}

\begin{proof}  The stronger statement follows since $f^{-1}(B_R(v))\subset X.$
So, it suffices to show $\as X\le n.$  Suppose $f$ is
$\lambda$-Lipschitz.  Given $\epsilon>0$ we construct an
$\epsilon$-Lipschitz, uniformly cobounded map $\psi:X\to K$ to a
uniform polyhedron of dimension $n.$

Let $c_{n-1}$ be the constant of uniformization from Proposition
\ref{Prop3}.  Take $\nu=\nu(\epsilon/4c_{n-1},n-1)$ and let
$r>\max\{\nu,8/\epsilon\}.$  Take a cover $\sW$ of $T$ by disjoint
sets so that the $\lambda r$-enlargement and the $2\lambda
r$-enlargement both have order $2.$

Since the $W\in \sW$ are uniformly bounded there is an $R>0$ so
that for each $W\in\sW$ there is a $v_W\in T$ so that $N_{2\lambda
r}(W)\subset B_R(v_W).$  Since $\as f^{-1}(B_R(v_W))\le n$
uniformly, $\as f^{-1}(N_{2\lambda r}(W)) \le n$ uniformly.

Consider a pair $W\neq W'$ for which $N_{\lambda r}(W)\cap
N_{\lambda r}(W')\neq\emptyset.$  By the finite union theorem,
$\as [f^{-1}(N_{\lambda r}(W))\cup f^{-1}(N_{\lambda r}(W'))]\le
n,$ so $\asI [f^{-1}(N_{\lambda r}(W))\cup f^{-1}(N_{\lambda
r}(W'))]\le n.$ Since $f^{-1}(W)$ and $f^{-1}(W')$ are
asymptotically disjoint, there is an asymptotic separator $A_e$
separating them in $f^{-1}(N_{\lambda r}(W))\cup f^{-1}(N_{\lambda
r}(W')),$ with $\asI A_e\le n-1.$ (Here, the subscript $e$ refers
to the edge $e=[W,W']$ in the nerve of $N_{\lambda r}(\sW).$)
Thus, $\as A_e\le n-1.$

As $N_r(A_e)$ is coarsely isometric to $A_e$ we have $\as
N_r(A_e)\le n-1.$  For each edge, let $\sV_e$ be a uniformly
bounded cover of $N_r(A_e)$ with multiplicity $\le n$ and with
Lebesgue number $L>r.$  For each $W$ cover $f^{-1}(N_{2\lambda
r}(W))$ by uniformly bounded sets with multiplicity $\le n+1$ and
with Lebesgue number greater than $\max\{b(\sV_e)\mid W\in e\},$
where $b(\sV_e)$ is an upper bound on the diameters of the sets in
$\sV_e.$  Since $X$ is assumed to have bounded geometry, this
maximum exists.

The conditions on the Lebesgue numbers along with the fact that
$N_r(A_e)\subset f^{-1}(N_{2\lambda r}(W))\cap f^{-1}(N_{2\lambda
r}(W'))$ guarantee that there exist simplicial maps
$g_W:Nerve(\sV_e)\to Nerve(\sU_W)$ and $g_{W'}:Nerve(\sV_e)\to
Nerve(\sU_{W'}).$  Take $M_{e,W}$ and $M_{e,W'}$ to be the uniform
mapping cylinders of the maps $g_W$ and $g_{W'},$ respectively.

As $r>8/\epsilon,$ we may apply Lemma \ref{L1} to $A_e\subset
\Gamma$ and the covers to obtain $\epsilon$-Lipschitz maps
$h_{e,W}: N_r(A_e)\to M_{e,W}$ and $h_{e,W'}: N_r(A_e)\to
M_{e,W'}$ to the uniform mapping cylinders.

For each $W\in\sW,$ construct a uniformly cobounded
$\epsilon$-Lipschitz map $\phi_W:f^{-1}(N_{2\lambda r}(W))\to
K_W,$ to the uniform $n$-dimensional simplicial complex $K_W$ by
taking the natural projection to the nerve of $\sV_W.$ Such a
mapping exists since $r>\nu,$ by Proposition \ref{nu}.

We note that the $N_r( A_e)$ are disjoint for distinct edges in
the nerve. Thus, for each $W\in \sW$ define
$\psi_W:\pi^{-1}(N_{2\lambda r}(W))\to K_W\cup_{W\in e}
M_{e,W}=L_W$ to the uniform complex $L_W,$ with mapping cylinders
attached as the union of the map $\phi_W$ restricted to
${f^{-1}(N_{2\lambda r}(W))\setminus\cup_{W\in e} N_r(A_e)}$ and
the restrictions of $h_{e,W}$ to ${N_r(A_e)\cap f^{-1}(N_{2\lambda
r}(W))},$ for all edges $e$ in $Nerve(N_{\lambda r}(\sW))$ which
contain $W$ as a vertex.

We construct $K$ by gluing together the $L_W.$  Clearly, the
dimension of $K$ is at most $n.$  The maps $\psi_W:\Gamma\to K$
agree on the common parts $A_e$ so they define a map $\psi:X\to
K.$ The map $\psi$ is $\epsilon$-Lipschitz by Proposition
\ref{Prop4}, and uniformly cobounded by Lemma \ref{L1}.
\end{proof}

What follows is a natural generalization of the combinatorial
structure of amalgamated free products studied in \cite{BD1}.

Let $X$ and $Y$ be pointed metric spaces.  Define a metric space
$X\hat{\ast}Y$ to be the metric space whose elements are
alternating words formed from the alphabet $X\setminus\{x_0\}
\sqcup Y\setminus\{y_0\}.$  Set $x_0=y_0=\tilde e.$  Define a norm
by the following rule: $\|z\|=0\iff z=\tilde e,$ and
$\|x_1y_1\cdots x_ry_r\|=\sum_id_X(x_i,x_0)+d_Y(y_i,y_0),$ where
we allow $x_1=x_0$ or $y_r=y_0.$  To define the metric, let $z,z'$
be words in $X\hat\ast Y.$ Write $z=uv$ and $z'=uv',$ so that $u$
is a common beginning, which we allow to be $\tilde e.$  Then,
$d(z,z')=\|v\|+\|v'\|.$  Observe that if $X$ and $Y$ are discrete
metric spaces with bounded geometry, then so is $X\hat\ast Y.$

We will need the following union theorems in the next proposition.
Both are taken from \cite{BD1}.

\begin{thm} {\em(Infinite Union Theorem)} Let $X_\alpha$ be a family of subsets of the metric
space $X$ satisfying the inequality $\as X_\alpha\le n$ uniformly.
Suppose that for every $r>0$ there exists a $Y_r\subset X$ so that
$\as Y_r\le n$ and the family $\{X_\alpha\setminus Y_r\}$ is
$r$-disjoint.  Then, $\as \cup_\alpha X_\alpha\le n.$
\end{thm}

By taking the family to consist of two sets $A$ and $B$ and taking
$Y_r=B$ for each $r$ we immediately obtain the following Finite
Union Theorem as a corollary.

\begin{thm} {\em(Finite Union Theorem)} Let $A$ and $B$ be subsets of a metric space $X.$  Then,
$\as A\cup B\le \max\{\as A,\as B\}.$
\end{thm}

\begin{Proposition} \label{XY} Let $(XY)^m$ denote the subset $XY\cdots
XY\subset X\hat\ast Y.$  Suppose that $\as X\le n$ and $\as Y\le
n.$  Then $\as (XY)^m\le n$ for all $m.$
\end{Proposition}

\begin{proof} Let $w\in X\hat\ast Y$ be a word.  Put $\ell(w)$ equal to the
length of $w,$ i.e. $k$ where $w=z_1\cdots z_k$ and the $z_i$
alternate, coming from $X$ and $Y.$ Put $P_k=\{w\mid \ell(w)=k\}.$
Denote by $P^X_k$ the set $\{w\in P_k\mid w_{\ell(w)}\in X\}.$
Similarly, put $P^Y_k=\{w\in P_k\mid w_{\ell(w)}\in Y\}.$  Since
$(XY)^m\subset \cup_{k=1}^{2m} P_k,$ by the finite union theorem,
it suffices to show that $\as P_k\le n$ for all $k.$

We proceed inductively.  If $k=1,$ then $P_k=X\cup Y$ so by the
finite union theorem $\as P_1\le n.$  Obviously, $P^X_{k+1}\subset
P_k^Y X$ and $P^Y_{k+1}\subset P^X_k Y.$  We show that $\as
P^X_{k+1}\le n.$  The other case is similar.  Put
$C_r=P^Y_kB_r^X(x_0),$ where $B^X_r(x_0)$ is the $r$-ball around
$x_0$ in $X.$  Then, $C_r\subset N^{X\hat\ast Y}_r(P_k),$ so $C_r$
is coarsely isometric to $P_k.$  Applying the inductive
hypothesis, we conclude that $\as C_r\le n.$

Next, consider the families $zX,$ where $z\in P^Y_{k}.$  Clearly
if $z\neq z',$ then $d(zx,z'x')>\|x\|+\|x'\|.$  Thus,
$\{zX\setminus C_r\}$ is an $r$ disjoint family.  Next, since for
every $z,$ $x\mapsto zx$ is an isometry in $X\hat\ast Y,$ the
families $zX$ are isometric.  Next, as $zX$ is coarsely isometric
to $X$ for all $z,$ we conclude that $\as zX\le n$ uniformly.  By
the infinite union theorem, we conclude that $\as P_{k+1}\le n.$
\end{proof}

Obviously the result of the previous theorem also holds for
subsets of the form $(YX)^m.$

There is a natural tree, $T,$ associated to $X\hat\ast Y.$ Define
the vertices of $T$ to be formal cosets $uX$ and $vY$ where $u$
and $v$ are words in $X\hat\ast Y.$  Connect the vertices $uX$ and
$vY$ by an edge if either $ux=v$ or $vy=u$ for some $x\in X,$ or
some $y\in Y.$

\begin{Proposition} As defined above, $T$ is a tree.
\end{Proposition}

\begin{proof}  Obviously $T$ is connected: given two vertices one
can find a path connecting them by starting at either of the root
vertices $x_0Y$ or $y_0X.$  Next if there were a circuit, say
$uX=ux_1y_1\cdots x_ry_r X,$ then this would mean that there exist
$x,x'\in X$ for which $d(x, x_1y_1\cdots x_ry_rx')=0.$  But $d(x,
x_1y_1\cdots x_ry_rx')\ge \sum_{i=2}^rd(x_0,x_i)+\sum_{i=1}^r
d(y_0,y_{i}).$ For this to be zero, we need all $x_i=\tilde e$ and
$y_i=\tilde e.$
\end{proof}

\begin{Theorem} \label{freeprod}
Let $X$ and $Y$ be discrete pointed metric spaces with bounded geometry, $\as
X= n,$ and $\as Y\le n,$ where $n>0.$  Then, $\as X\hat\ast Y=n.$
\end{Theorem}

\begin{proof}  By Theorem \ref{totrees} we must find a Lipschitz
map to a tree, show that bounded disjoint sets in the tree lift to
asymptotically disjoint sets in $X\hat\ast Y$ and that for every
$R,$ $\as f^{-1}(B_R(v))\le n$ for all $v\in T.$

Take $T$ as above, and define $f:X\hat\ast Y\to T$ by $f(u)=uX.$
Let $u,v\in X\hat\ast Y$ with common part $w.$  Then, $u=wu'$ and
$v=wv'.$  Then, $d(u,v)=\|u'\|+\|v'\|\ge \ell(u')+\ell(v')$ where
$\ell(t)$ is the length of $t.$  Since
\[d(f(u),f(v))=d(uX,vX)=d(u'X,v'X)=\ell(u')+\ell(v'),\]
we conclude that $f$ is $1$-Lipschitz.

Next, observe that if $uX$ and $vX$ are distinct vertices of $T,$
then $d(ux,vx')\ge \|x\|+\|x'\|.$  Let $W$ and $W'$ be disjoint
bounded subsets of $T.$  Then, the sets $\{ux_0\colon uX\in W\}$
and $\{vx_0\colon vX\in W'\}$ are bounded.  Let $t_0X\in T$ be
given and take $r$ so large that $B_{r/2}(t_0X)\subset X$ contains
$\{ux_0\colon uX\in W\}\cup \{vx_0\colon vX\in W'\}.$ Then,
$d(f^{-1}(W)\setminus B_r(t_0X),f^{-1}(W')\setminus B_r(t_0X))\ge
r.$  Thus, $f^{-1}(W)$ and $f^{-1}(W')$ are asymptotically
disjoint.

Finally, it is easy to see that $f^{-1}(B_R(vX))\subset vX(YX)^R,$
and so, by Proposition \ref{XY}, $\as B_R(vX)\le n.$
\end{proof}

\begin{Theorem} \label{amalgam}
Let $A$ and $B$ be finitely generated groups with finite
asymptotic dimension.  Let $C$ be a common subgroup.  Then $\as
A\ast_C B\le \as C+\max\{\as A/C, \as B/C, 1\}.$
\end{Theorem}

{\bf Remark:} This estimate is not always an improvement over the
previously known estimate $\as A\ast_C B\le 1+\max\{\as A, \as
B\},$ (see \cite{BD2}) since there is no way to give an upper
bound on $\as A/C$ in terms of $\as A$ and $\as C.$  In
particular, Thompson's group $F$ is a two generator group, and
hence is a quotient of $\mathbb{F}_2;$ but $\as F=\infty,$ (as it
contains a copy of $\mathbb{Z}^n$ for each $n$) whereas $\as
\mathbb{F}_2=1.$

\begin{proof}
It is well-known that every element $x\in A\ast_C B$ admits a
unique normal presentation $c\bar x_1\bar x_2\cdots \bar x_k$
where $c\in C$ and $\bar x_i=Cx_i$ are non-trivial alternating
cosets of $C$ in $A$ or $B$ and $x=cx_1\cdots x_k.$  Given a
metric on $A$ and $B$ we define a metric on $C\backslash A$ and
$C\backslash B$ by $\bar d(Cx,Cy)=d_{A\ast_C B}(x,Cy).$  Thus, we
can consider the metric space $(C\backslash A)\hat\ast
(C\backslash B)$ where the common point is $\tilde e.$

Define a map $\phi:A\ast_C B\to (C\backslash A)\hat\ast
(C\backslash B)$ by defining $\phi(e)=\tilde e$ and
$\phi(x)=x_1\cdots x_k$ where $x=cx_1\cdots x_k$ is the normal
presentation of $x.$ We claim that $\phi$ is $1$-Lipschitz.  Since
$A\ast_C B$ is a discrete geodesic metric space, it suffices to
check the Lipschitz condition on pairs of the form $(x,xs),$ where
$s$ is in the generating set $S.$  The normal presentation of $xs$
will be either $c\bar x_1\cdots \overline{x_ks},$ or $c\bar
x_1\cdots \bar x_k\bar s.$  In the first case,
$d(\phi(x),\phi(xs))=\bar d(Cx_k,Cx_k s)\le d(x_k,x_k s)=1.$  In
the second, $d(\phi(x),\phi(xs))=\bar d(C,Cs)\le 1.$

By Theorem \ref{freeprod}, $\as (C\backslash A)\hat\ast
(C\backslash B)\le \max\{\as (C\backslash A), \as (C\backslash
B),1\}.$

Consider $\phi^{-1}(B_{2R}(\tilde e)).$  First, observe that
$B_{2R}(\tilde e)$ consists of alternating words $x_1x_2\cdots
x_k$ where the $x_i$ alternate between $C\backslash A$ and
$C\backslash B,$ and $\|x_1\cdots x_k\|\le 2R.$ Thus,
$\phi^{-1}(B_{2R}(\tilde e))\subset \cup_{\|w\|\le 2R} Cw.$ Since
$(C\backslash A)\hat\ast(C\backslash A)$ has bounded geometry,
this is a finite union. Applying the finite union theorem we see,
$\as \phi^{-1}(B_{2R}(\tilde e))\le \as \cup_w Cw\le \max\{\as
Cw\colon \|w\|\le 2R\}.$  But, since $Cw$ is coarsely isometric to
$C,$ we obtain $\as \phi^{-1}(B_{2R}(\tilde e))\le \as C.$

To get the inequality $\as\phi^{-1}(B_R(x))\le \as C$ uniformly,
we appeal to Proposition 1 of \cite{BD1}, which says that $\as
F_\alpha\le n$ uniformly if there exist $1$-Lipschitz injective
maps $f_\alpha:F_\alpha\to X$ to a metric spaces with bounded
geometry and $\as X\le n.$

For each $x\in (C\backslash A)\hat\ast(C\backslash B),$ write
$x=\omega x'$ where either $\|x'\|=R$ or else $\omega=\tilde e.$
Suppose $y\in B_R(x).$  Then, let $z$ be the common part of $x$
and $y,$ so that $x=zx^{\prime\prime}$ and $y=zy^{\prime\prime}.$
Then, $d(x'',y'')\le R,$ and since $x''$ and $y''$ have no common
beginning, we conclude that $\|x''\|\le R.$  Thus, $z=\omega z'.$
Hence, $y=\omega z'y'',$ where $\|z'\|\le R$ and $\|y''\|\le R.$
We conclude that $y\in B_{2R}(\omega).$

Let $f_x:\phi^{-1}(B_R(x))\to \phi^{-1}(B_{2R}(\tilde e))$ be
defined by $y\mapsto y''$ where $y''$ is the word $y$ with the
beginning part $\omega$ removed from it.  Then $f_x$ is an
isometry into $\phi^{-1}(B_{2R}(\tilde e)),$ which is a bounded
geometry space with $\as \phi^{-1}(B_{2R}(\tilde e))\le \as C.$

Thus, by the Hurewicz theorem, we have the desired estimate.
\end{proof}

\section{An Extension Theorem for Asdim}

Although it was known for some time (see \cite{BD1}) that
extensions of groups with finite asymptotic dimension had finite
asymptotic dimension, the Hurewicz-type theorem for group actions,
Theorem \ref{main}, allows us to give a sharp upper bound estimate
for the dimension.

\begin{Theorem}  \label{extension} Let $\phi:G\to H$ be a surjective homomorphism of
a finitely generated group with kernel $K.$  Suppose that $\as
H\le n$ and $\as K\le k.$  Then, $\as G\le n+k.$
\end{Theorem}

\begin{proof}  Let $S$ be a finite generating set for $G,$ and
take the set $\phi(S)$ as a generating set for $H.$  We consider
$G$ and $H$ in the left-invariant word metric.  The group $G$ acts
on $H$ by isometries according to the rule $g.h=\phi(g)h.$

We claim that $W_R(e)=N_R(K),$ where $e$ is the identity element.
Indeed, if $d_S(g,K)\le R,$ then $d_{\phi(S)}(\phi(g),e)\le R.$ On
the other hand, if $g\in W_R(e),$ then $\|\phi(g)\|_{\phi(S)}\le
R.$ Let $\phi(g)=\phi(s_{i_1})\cdots \phi(s_{i_k}),$ where
$s_{i_j}\in S,$ and $k\le R.$ Then, $g s_{i_k}^{-1}\cdots
s_{i_1}^{-1}\in K,$ and $d_S(g,g s_{i_k}^{-1}\cdots
s_{i_1}^{-1})=k\le R.$

Since $N_R(K)$ is coarsely equivalent to $K,$ $\as N_R(K)\le k,$
and the result follows from the theorem.
\end{proof}

Using the extension theorem we can prove the following form of
Theorem \ref{amalgam} very easily.

\begin{Proposition} Let $C\triangleleft A$ and $C\triangleleft B$
where $A$ and $B$ are finitely generated groups with finite $\as.$
Then $\as A\ast_C B\le \as C+\max\{\as C\backslash A, \as
C\backslash B, 1\}.$
\end{Proposition}

\begin{proof}  There is a natural surjection of groups $A\ast_C B\to(C\backslash A)\ast
(C\backslash B)$ with kernel $C.$  By the extension theorem, $\as
A\ast_C B\le \as C+\as(C\backslash A)\ast (C\backslash B).$
Applying the formula for the $\as$ of a free product from
\cite{BDK}, we get $\as A\ast_C B\le \as C+\max\{\as C\backslash
A, \as C\backslash B, 1\}.$
\end{proof}

Recall that a group $G$ is called {\em polycyclic} if there exists
a sequence of subgroups $\{1\}=G_0\subset G_1\subset\cdots\subset
G_n=G$ such that each $G_i\lhd G_{i+1}$ and $G_{i+1}/G_i$ is
cyclic.  The {\em Hirsch length} of $G$, denoted $h(G)$ is the
number of factors $G_{i+1}/G_{i}$ isomorphic to $\mathbb{Z}.$

\begin{Theorem} Let $\Gamma$ be a finitely generated polycyclic
group.  Then $\as \Gamma\le h(\Gamma).$
\end{Theorem}

\begin{proof}  Denote the
sequence of subgroups satisfying the polycyclic condition by:
$\{1\}=\Gamma_0\subset \Gamma_1\subset\cdots\subset
\Gamma_n=\Gamma.$ Then, by Theorem \ref{extension}, we have

\[\begin{array}{rcl}
\as\Gamma & \le&
\as\Gamma_n/\Gamma_{n-1}+\Gamma_{n-1}\\
&\le& \as\Gamma_n/\Gamma_{n-1}+\as\Gamma_{n-1}/\Gamma_{n-2}+\as
\Gamma_{n-2}\\
&\vdots&\\
&\le &\as\Gamma_n/\Gamma_{n-1}+\cdots+\as\Gamma_{1}/\Gamma_{0}+\as
\Gamma_{0}.
\end{array}\]
Since $\as \Gamma_{i+1}/\Gamma_i$ is only positive if
$\Gamma_{i+1}/\Gamma_i$ is isomorphic to $\mathbb{Z},$ and since
in this case, $\as \Gamma_{i+1}/\Gamma_i=1,$ we conclude $\as
\Gamma\le h(\Gamma).$
\end{proof}

Since every finitely generated nilpotent group is polycyclic we
immediately obtain the following result.

\begin{Corollary} \label{6} Let $\Gamma$ be a finitely generated nilpotent
group.  Then $\as\Gamma\le h(\Gamma).$
\end{Corollary}

Corollary \ref{6} can be extended to nilpotent Lie groups $N$ if
one defines the Hirsch length $h(N)$ as the sum of the number of
factors in $\Gamma_{i+1}/\Gamma_i$ isomorphic to $\mathbb{R}$ for
the central series $\{\Gamma_i\}$ of $N$. We take an equivariant
metric on $N$ and on the quotients. Then the projection
$\Gamma_{i+1}\to \Gamma_{i+1}/\Gamma_i$ is 1-Lipschitz and
$\Gamma_{i+1}/\Gamma_i$ is coarsely isomorphic to
$\mathbb{R}^{n_i}$. Then we have

\begin{Corollary}  \label{Lie} Let $N$ be a nilpotent Lie group endowed with an equivariant
metric. Then $\as N\le h(N).$
\end{Corollary}

Since $h(N)=\dim N$ for simply connected $N$, we obtain

\begin{Corollary} \cite[Theorem 3.5]{CG}  \label{cg} For
a simply connected nilpotent Lie group $N$ endowed with an
equivariant metric $\as N\le \dim N.$
\end{Corollary}

Actually in view of \cite[Corollary 1.F1]{Gr93} the inequalities
in Corollaries \ref{Lie} and \ref{cg} are equalities.

Corollary \ref{cg} is the main step in the proof of the following

\begin{Theorem}\cite{CG} For a connected Lie group $G$ and its maximal compact subgroup $K$
there is a formula $\as G/K=\dim G/K$ where $G/K$ is endowed with
a G-invariant metric.
\end{Theorem}
This theorem in particular allows to compute asymptotic dimension
of the hyperbolic space $\mathbb{H}^n=n$.

\begin{Corollary} $\as \mathbb{H}^n=n$.
\end{Corollary}

\begin{proof}
Take $G=O(n,1)_+$ and $K=O(n)$.
\end{proof}
This computation can be generalized in spirit of \cite{Ro3}.

Let $(X,d)$ be a metric space. By $\sH(X)$ we denote the space of
balls in $X$ endowed with the following metric
$$
\rho(B_t(x),B_s(y))=2\ln(\frac{d(x,y)+\max\{t,s\}}{\sqrt{ts}}).
$$

We note that $\sH(\mathbb{R}^n)$ is coarsely equivalent to
$\mathbb{H}^{n+1}$, \cite[Example 2.60]{Ro3}.

We recall that a  metric space $X$ with $\as X\le n$ is said to
satisfy the {\em Higson property} \cite{DrZ} if there exists $C>0$
such that for every $D>0$ there exists a cover $\sU$ of $X$ with
$\mesh(\sU)<CD$ and such that $\sU=\sU^0\cup\dots\cup\sU^n$, where
$\sU^0,\dots,\sU^n$ are $D$-disjoint. In \cite{Ro3} $X$ satisfying
this condition are said to have asymptotic dimension $\le n$ {\em
of linear type}. It is shown in \cite{DrZ} that every metric space
of bounded geometry with $\as X\le n$ admits a coarsely equivalent
metric with the Higson property. Unfortunately the coarse type of
$\sH(X)$ depends on a metric on $X$ not only the coarse class of
metrics.

\begin{Theorem} Suppose that the metric space $(X,d)$
possesses the Higson property.
Then $\as\sH(X)=\as X+1$.
\end{Theorem}
\begin{proof}
Consider the projection $\pi:\sH(X)\to\mathbb{R}$ defined by
$\pi(B_t(x))=\ln t$ and apply Theorem \ref{Hurewicz} to it (see
\cite[Corollary 9.21]{Ro3}).
\end{proof}


\begin{thebibliography}{10}

\bibitem{Be04}
G.~Bell, \emph{Asymptotic properties of groups acting on
complexes}, to appear
  in Proc. Amer. Math. Soc., 2003.

\bibitem{BD1}
G.~Bell and A.~Dranishnikov, \emph{On asymptotic dimension of
groups}, Algebr.
  Geom. Topol. \textbf{1} (2001), 57--71.

\bibitem{BD2}
\bysame, \emph{On asymptotic dimension of groups acting on trees},
Geom.
  Dedicata \textbf{103} (2004), 89--101.

\bibitem{BDK}
G.~Bell, A.~Dranishnikov, and J.~Keesling, \emph{On a formula for
the
  asymptotic dimension of free products}, Submitted, 2004.

\bibitem{BH}
M.~Bridson and A.~Haefliger, \emph{Metric spaces of non-positive
curvature},
  Springer, 1999.

\bibitem{CG}
G.~Carlsson and B.~Goldfarb, \emph{On homological coherence of
discrete
  groups}, J. Algebra \textbf{276} (2004), 502--514.

\bibitem{Dr00}
A.~Dranishnikov, \emph{Asymptotic topology}, Russian Math. Surveys
\textbf{55}
  (2000), no.~6, 71--116.

\bibitem{Dr01}
\bysame, \emph{On asymptotic inductive dimension}, JP Jour.
Geometry \&
  Topology \textbf{1} (2001), no.~3, 239--247.

\bibitem{DJ}
A.~Dranishnikov and T.~Januszkiewicz, \emph{Every coxeter group
acts amenably
  on a compact space}, Topology Proc. \textbf{24} (1999), 135--141.

\bibitem{DrZ}
A.~Dranishnikov and M.~Zarichnyi, \emph{Universal spaces for
asymptotic
  dimension}, Topology Appl. \textbf{140} (2004), 203--225.

\bibitem{En95}
R.~Engelking, \emph{Theory of dimensions finite and infinite},
Sigma Series in
  Pure Mathematics, vol.~10, Heldermann Verlag, 1995.

\bibitem{Gr93}
M.~Gromov, \emph{Asymptotic invariants of infinite groups},
Geometric Group
  Theory, London Math. Soc. Lecture Note Ser. (G.~Niblo and M.~Roller, eds.),
  no. 182, 1993.

\bibitem{Ji}
L.~Ji, \emph{Asymptotic dimension of arithmetic groups.}, preprint
(2003).

\bibitem{Ro3}
J.~Roe, \emph{Lectures on coarse geometry}, University Lecture
Series, vol.~31,
  AMS, 2003.

\bibitem{Yu1}
G.~Yu, \emph{The {N}ovikov conjecture for groups with finite
asymptotic
  dimension}, Annals of Mathematics \textbf{147} (1998), no.~2, 325--355.

\end{thebibliography}

\end{document}